\documentclass{amsart}

\usepackage{amssymb}
\usepackage{amsmath, amsfonts, enumerate}
\usepackage{color}
\usepackage{amsthm}
\usepackage{hyperref}
\usepackage{mathrsfs}
\usepackage{txfonts}

\begin{document}

\newtheorem{theorem}{Theorem}    %[section]
\newtheorem{proposition}[theorem]{Proposition}
\newtheorem{conjecture}[theorem]{Conjecture}
\newtheorem{corollary}[theorem]{Corollary}
\newtheorem{lemma}[theorem]{Lemma}
\newtheorem{sublemma}[theorem]{Sublemma}
\newtheorem{observation}[theorem]{Observation}
\newtheorem{remark}[theorem]{Remark}
\newtheorem{definition}[theorem]{Definition}
\newtheorem{notation}[theorem]{Notation}
\newtheorem{question}[theorem]{Question}
\newtheorem{questions}[theorem]{Questions}
\newtheorem{example}[theorem]{Example}
\newtheorem{problem}[theorem]{Problem}
\newtheorem{exercise}[theorem]{Exercise}

\numberwithin{theorem}{section}
\numberwithin{theorem}{section}
\numberwithin{equation}{section}

\newcommand{\ri}{{\rm{i}}}
\newcommand{\inner}[1]{ \langle  #1 \rangle}
\newcommand{\Norm}[1]{ \left\|  #1 \right\| }
\newcommand{\set}[1]{ \left\{ #1 \right\} }
\newcommand{\modulo}[2]{[#1]_{#2}}
\newcommand{\bra}[1]{\Bigg(  #1 \Bigg)}
\newcommand{\cir}[1]{\left(  #1 \right)}
\newcommand{\wh}[1]{\widehat{#1}}
\newcommand{\wt}[1]{\widetilde{#1}}
\newcommand{\cl}[1]{{\rm  #1 }}
\newcommand{\lp}[1]{L_{#1}(\ted)}
\newcommand{\np}[1]{L_{#1}(\mathcal{M})}
\newcommand{\nlp}[1]{L_{#1}(\RR^{d};\A_{\theta}^{d})}

\def\one{\mathbf 1}
\def\sn{{{S}^{n-1}}}
\def\qd{\,{\mathchar'26\mkern-12mu d}}
\def\p{\partial}
\def\rp{ ^{-1} }
\def\supp{{\rm supp}}
\def\ted{\mathbb{T}_{\theta}^{d}}
\def\td{\mathbb{T}^{d}}
\def\zd{\mathbb{Z}^{d}}
\def\rd{\mathbb{R}^{d}}
\def\cci{C_{c}^{\infty}(\rd)}
\def\e{\mathbf e}
\def\Bl{\B(\lp{2})}
\def\BMO{{\rm BMO}}
\def\pair{\RR^{d};\A_{\theta}^{d}}
\def\atd{\C(\ted)}
\def\bsv{\dot{B}_{\infty}\rp(\pair)}
\def\proj{{\rm Proj}}
\def\lo{L_{0}(\ted)}
\def\dpr{\prec\prec}

\def\MM{\mathbb{M}}
\def\NN{\mathbb{N}}
\def\QQ{{\mathbb{Q}}}
\def\SS{\mathbb{S}}
\def\RR{{\mathbb R}}
\def\TT{{\mathbb T}}
\def\HH{{\mathbb H}}
\def\ZZ{{\mathbb Z}}
\def\CC{{\mathbb C}}
\def\EE{\mathbb E}
\def\GG{\mathbb G}

\def\VV{\mathbf V}

\def\H{{\mathcal H}}
\def\L{\mathcal{L}}
\def\S{\mathcal{S}}
\def\F{{\mathcal F}}
\def\G{{\mathcal G}}
\def\M{{\mathcal M}}
\def\B{{\mathcal B}}
\def\C{{\mathcal C}}
\def\T{{\mathcal T}}
\def\I{{\mathcal I}}
\def\E{{\mathcal E}}
\def\V{{\mathcal V}}
\def\W{{\mathcal W}}
\def\U{{\mathcal U}}
\def\A{{\mathcal A}}
\def\R{{\mathcal R}}
\def\O{{\mathcal O}}
\def\D{{\mathcal D}}
\def\N{{\mathcal N}}
\def\P{{\mathcal P}}
\def\K{{\mathcal K}}

\title[Gagliardo--Nirenberg interpolation inequality]{Gagliardo--Nirenberg interpolation inequality for symmetric spaces on Noncommutative torus}

%11M55,


\author{Fedor SUKOCHEV}
\address{School of Mathematics and Statistics, UNSW, Kensington, NSW 2052, Australia}
\email{f.sukochev@unsw.edu.au}

\author{Fulin YANG}
\address{School of Mathematics, Harbin Institute of Technology, Harbin 150001, China}
\email{fulinyoung@hit.edu.cn}

\author{Dmitriy ZANIN}
\address{School of Mathematics and Statistics, UNSW, Kensington, NSW 2052, Australia}
\email{d.zanin@unsw.edu.au}

%\markboth{F. A. Sukochev, F. L. Yang and D. Zanin}
%{Gagliardo-Nirenberg inequality}

\date{}
\maketitle

\begin{abstract}
Let $E(\ted),F(\ted)$ be two symmetric operator spaces on noncommutative torus $\ted$ corresponding to symmetric function spaces $E,F$ on $(0,1)$.
We obtain the Gagliardo--Nirenberg interpolation inequality  with respect to $\ted$: if  $G=E^{1-\frac{l}{k}}F^{\frac{l}{k}}$ with $ 0\leq l\leq k$ and if the Ces\`{a}ro operator is bounded on $E$ and $F$, then
\begin{align*}
\|\nabla^lx\|_{G(\ted)}\leq 2^{3\cdot 2^{k-2}-2}(k+1)^d\|C\|_{E\to E}^{1-\frac{l}{k}}\|C\|_{F\to F}^{\frac{l}{k}}\|x\|_{E(\ted)}^{1-\frac{l}{k}}\|\nabla^kx\|_{F(\ted)}^{\frac{l}{k}},\; x\in W^{k,1}(\ted),
\end{align*}
where $W^{k,1}(\ted)$ is the Sobolev space on $\ted$ of order $k\in\NN$.
%Especially, the assumption of \cite[Theorem 2.1]{LRS2023} on Boyd indices of the given symmetric function spaces $E$ and $F$ can be dispensed of. 
Our method is different from the previous settings, which is of interest in its own right.
\end{abstract}

%%%%%%%%%%%%%%%%%%%%%%%%%%%%%%%%%%%%%%%%%%%%%%%%%%%%%%%%%%%%%%%%%%%%%%%
\section{Introduction}
The Gagliardo--Nirenberg inequality has wide applications in various areas, such as the solvability of partial differential equation \cite{BDGO1997,PD2002,PP2004}, the study of the chain rule for Sobolev spaces \cite{MRS2021} and the multiplier theory \cite{GHS2020,S2019,XXY2018}.
At the beginning of the 20th century, Landau \cite{L1914} obtained the following inequality
\begin{align*}
\Norm{u'}_{L_{\infty}(\RR)}\lesssim\Norm{u''}_{L_{\infty}(\RR)}^{\frac{1}{2}}\Norm{u}_{L_{\infty}(\RR)}^{\frac{1}{2}},\quad u\in W^{2,\infty}(\RR).
\end{align*}
In the sequel, we write $A\lesssim B$ if $A\leq cB$ for some constant $c>0$ which does not depend on $A$ and $B$, and $W^{k,p}(\Omega)$ $(k\in\NN,1\leq p\leq\infty)$ denotes the Sobolev space on $\Omega$. 
Later on, Kolmogoroff \cite{K1949} proved the result for higher-order derivatives. In 1958,  Nash \cite{N1958} showed that
\begin{align*}
\Norm{u'}_{L_{4}(\Omega)}\lesssim\Norm{\nabla u}_{L_{2}(\Omega)}^{\frac{1}{2}}\Norm{u}_{L_{2}(\Omega)}^{\frac{1}{2}},\quad u\in W^{1,2}(\Omega),~\Omega=\RR^{2}~\mbox{or}~\RR^{3}.
\end{align*}
Stein \cite{S1957} generalized the results due to Landau and Kolmogoroff to the setting of $L_{p}(\RR)$ by showing that
\begin{align*}
\Norm{u'}_{L_{2p}(\RR)}\lesssim\Norm{u''}_{L_{p}(\RR)}^{\frac{1}{2}}\Norm{u}_{L_{p}(\RR)}^{\frac{1}{2}}
,\quad u\in W^{2,p}(\RR),~ 1\leq p<\infty.
\end{align*} 
In 1959, Gagliardo \cite{G1958} and Nirenberg \cite{N1959} independently established the following result
\begin{align*}
\Norm{\nabla^{l}u}_{L_{p}(\RR^{d})}
\lesssim\Norm{\nabla^{k}u}_{L_{r}(\RR^{d})}^{\frac{l}{k}}\Norm{u}_{L_{q}(\RR^{d})}^{1-\frac{l}{k}}
\quad\mbox{with}\quad 0\leq l<k,~\frac{k}{p}=\frac{l}{r}+\frac{k-l}{q}.
\end{align*}
for $1\leq r<\infty$, $1\leq q\leq\infty$ and $u\in W^{k,r}(\RR^{d})\bigcap L_{q}(\RR^{d})$. This inequality is now called the Gagliardo--Nirenberg interpolation inequality.
In \cite{XXY2018}, a noncommutative version of the Gagliardo--Nirenberg interpolation inequality was established in $L_{p}$-space $(1<p<\infty)$ of the noncommutative torus $\ted$ (see Subsection \ref{noncommutative torus} for the definition of $\ted$).
Recently, Le\'{s}nik et al. \cite{LRS2023} further generalized the Gagliardo--Nirenberg interpolation inequality, by replacing the function spaces with symmetric function spaces $E(\RR^{d})$ and $F(\RR^{d})$, i.e.
\begin{align}\label{GN with Boyd indices}
\Norm{\nabla^{l}u}_{E(\RR^{d})^{\frac{l}{k}}F(\RR^{d})^{1-\frac{l}{k}}}
\lesssim\Norm{\nabla^{k}u}_{E(\RR^{d})}^{\frac{l}{k}}\Norm{u}_{F(\RR^{d})}^{1-\frac{l}{k}}
\quad\mbox{with}\quad 0\leq l<k,u\in W^{k,1}_{loc}(\RR^{d})\cap F(\RR^{d}),
\end{align}
where $E(\RR^{d})^{\frac{l}{k}}F(\RR^{d})^{1-\frac{l}{k}}$ is defined by Calder\'{o}n-Lozanovskii construction (see e.g \cite{DPS2024}, also Subsection \ref{Calderon-Lozanovskii space} below) and the upper Boyd indices (see e.g. \cite{PKJF2013}) of $E(\RR^{d})$, $F(\RR^{d})$ are required to be strictly smaller than 1.
The main tools used in \cite{LRS2023} are Hardy-Littlewood maximal functions on $\RR^{d}$ and the upper Boyd index assumptions guarantee that maximal functions are bounded on  $E(\RR^{d})$ and $F(\RR^{d})$. 
%Le\'{s}nik et al. \cite[Question 4.1]{LRS2023} asked the following question: {\it does \eqref{GN with Boyd indices} hold without the assumptions on Boyd indices?} %In this paper, we give a positive answer to this question.

%For our purpose, we use the notation of the Ces\`{a}ro operator $C:L_{1}(0,1)\rightarrow L_{0}(0,1)$ which is defined as
%\begin{align*}	C(f)(t)=\frac{1}{t}\int_{0}^{t}f(s)ds,\quad f\in L_{1}(0,1).
%\end{align*}

Our main result is the following Theorem \ref{GN-symmetric}: Gagliardo--Nirenberg interpolation inequality for the symmetric operator spaces on noncommutative torus $\ted$. 
%One of the main purpose of this paper is to answer \cite[Question 4.1]{LRS2023} in the affirmative. Indeed, our main result, Theorem \ref{GN-symmetric} below, resolves this question even in the noncommutative framework. 
The Ces\`{a}ro  operator plays an important role in the proof of our main result. Recall that the Ces\`{a}ro operator $C$ on the space $L^{1}_{loc}(\RR_{+})$ of the locally integrable functions on $\RR_{+}$ is defined by
\begin{align*}
C(\varphi)(t)=\frac{1}{t}\int_{0}^{t}\varphi(s)ds,\quad t>0,\quad\varphi\in L^{1}_{loc}(\RR_{+}).
\end{align*}
Let $\p_{j}$ be the $j$-th partial derivative on $\ted$ and let $\nabla^{\alpha}=\p_{1}^{\alpha_{1}}\cdots\p_{d}^{\alpha_{d}}$ for $\alpha=(\alpha_{1},\cdots,\alpha_{d})\in\NN^{d}$ (see Subsection \ref{noncommutative torus} for $\ted$). Denoting $|\alpha|_1=\alpha_{1}+\cdots+\alpha_{d}$, we use the notation
$$\|\nabla^kx\|_{G(\ted)}=\max_{|\alpha|_1=k}\|\nabla^{\alpha}x\|_{G(\ted)},$$
for the symmetric operator space $G(\ted)$ on $\ted$ corresponding to the symmetric function space $G$ on $(0,1)$.
Suppose that $E^{\epsilon}F^{1-\epsilon}$ is defined by the Calder\'{o}n--Lozanovskii construction for symmetric function spaces $E,F$ on $(0,1)$ and $0<\epsilon<1$ (see Subsection \ref{Calderon-Lozanovskii space}).

\begin{theorem}\label{GN-symmetric}
Let $0\leq l\leq k$. Suppose that $E,F$ are two symmetric function spaces on $(0,1)$ and set $G=E^{\frac{l}{k}}F^{1-\frac{l}{k}}.$ If the Ces\`{a}ro operator $C$ is bounded on $E$ and on $F,$ then
$$\|\nabla^lx\|_{G(\ted)}\leq 2^{3\cdot 2^{k-2}-2}(k+1)^d\|C\|_{E\to E}^{1-\frac{l}{k}}\|C\|_{F\to F}^{\frac{l}{k}}\|x\|_{E(\ted)}^{1-\frac{l}{k}}\|\nabla^kx\|_{F(\ted)}^{\frac{l}{k}},\quad x\in W^{k,1}(\ted),$$
where $W^{k,1}(\ted)$ is the Sobolev space on $\ted$ of order $k\in\NN$.
\end{theorem}
%Due to lack of suitable notion of the maximal function, our approach, originally designed for $\ted$, allows us to provide a new proof for the purely commutative cases without using maximal functions. 
%It is worth noting that, in Theorem \ref{GN-symmetric}, we provide a constant in the inequality with respect to the symmetric function spaces and the Ces\`{a}ro operator. It may be helpful in partial differential equation.
\begin{remark}
Let $0\leq l\leq k$ and $x\in W^{k,1}(\ted)$. Abandoning the use of Ces\`{a}ro operator in the proof of Theorem \ref{GN-symmetric}, the following interpolation inequality also holds
$$\|\nabla^lx\|_{L_{1}(\ted)}\leq 2^{3\cdot 2^{k-2}-2}(k+1)^d
\|x\|_{L_{1}(\ted)}^{1-\frac{l}{k}}\|\nabla^kx\|_{L_{1}(\ted)}^{\frac{l}{k}}.$$
\end{remark}

Recall that the approach used in \cite{LRS2023} strongly relies on the Hardy-Littlewood maximal function, which is not applicable in the noncommutative setting. Our method is completely different from that used in \cite{LRS2023}, which is of interest in its own right.

Let us explain the idea used in our paper in detail. To overcome the difficulty in establishing the Gagliardo--Nirenberg inequality for symmetric spaces on $\ted$, it is worth to revisit the methods used in the setting of $\RR^{d}$. On $\RR^{d}$, the pointwise version of the Gagliardo--Nirenberg inequality in terms of the maximal functions were obtained by (see e.g. \cite{K1994,MS1999,LRS2023}) i.e.
\begin{align}\label{nabla maximal funxtion}
	|\nabla^{l}u|\lesssim M(\nabla^{k}u)^{\frac{l}{k}}M(u)^{1-\frac{l}{k}},\quad u\in W^{k,1}_{loc}(\RR^{d}),
\end{align}
where $M$ is the Hardy-Littlewood maximal function over balls defined by
$\displaystyle M(f)(\xi)=\sup_{Q\ni\xi}\frac{1}{|Q|}\int_{Q}|f(y)|dy$ for $\xi\in\RR^{d}$, ball $Q\subset\RR^{d}$ and $f\in L_{loc}^{1}(\RR^{d})$.
%see e.g. \cite{G2014-c}.
By the {\it Riesz-Herz equivalence} \cite[Theorem 3.8]{BS1988}, inequality \eqref{nabla maximal funxtion} can be rewritten in the following form
\begin{align}\label{intermidiate ineq}
	C\mu(\nabla^{l}u)\lesssim \left(C\mu(\nabla^{k}u)\right)^{\frac{l}{k}} \left(C\mu(u)\right)^{1-\frac{l}{k}},\quad u\in W^{k,1}_{loc}(\RR^{d}),
\end{align}
where $\mu$ is the non-increasing rearrangement function (or singular valued function) (see Subsection \ref{nocommutative symmetric spaces}). This inequality does not involve maximal functions. In the setting of $\ted$, we establish a noncommutative version of \eqref{intermidiate ineq} for the case when $l=1$ and $k=2$. 
The following equality
$$iA=e^{iA}-{\rm id}+\int_{0}^{1}\int_{0}^{t} A^2e^{isA}dsdt$$
(here, $A$ is a self-adjoint densely defined closed operator on a Hilbert space) is crucial to establish, for $x\in W^{k,1}(\ted)$,
\begin{align*}%\label{important eq}	
C\mu(D_{j}x)\leq2 \left(C\mu(D_{j}^{2}x)\right)^{\frac{1}{2}} \left(C\mu(x)\right)^{\frac{1}{2}}
\quad\mbox{with}\quad D_{j}=i\p_{j},~j=1,\cdots,d,
\end{align*}
which is the key to prove our main result.

In particular, when $\theta=0$, $\ted$ is the commutative von Neumann algebra $L_{\infty}(\TT^{d})$. By Theorem \ref{GN-symmetric}, we obtain the Gagliardo--Nirenberg interpolation inequality for the symmetric function spaces on $\TT^{d}$. 
\begin{corollary}\label{commutative torus}
Let $E(\TT^{d})$ and $F(\TT^{d})$ be
the symmetric function spaces on $\TT^{d}$. If the Ces\`{a}ro operator $C$ is bounded on $E(\TT^{d})$ and $F(\TT^{d})$, then
\begin{align*}
\Norm{\nabla^{l}u}_{E(\TT^{d})^{\frac{l}{k}}F(\TT^{d})^{1-\frac{l}{k}}}
\lesssim\Norm{\nabla^{k}u}_{E(\TT^{d})}^{\frac{l}{k}}\Norm{u}_{F(\TT^{d})}^{1-\frac{l}{k}}
\quad\mbox{with}\quad u\in W^{k,1}(\TT^{d})\cap F(\TT^{d})
,~ 0\leq l<k.
\end{align*}
\end{corollary}
%This Corollary gives the Gagliardo--Nirenberg interpolation inequality for the symmetric function spaces on $\TT^{d}$ without any assumptions on Boyd indices.

\begin{remark}\label{classical RRd}
If we replace $\mathbb{T}^d$ with $\mathbb{R}^d$, a verbatim repeating of our argument in the proof still works on $\RR^{d}$. Thus we can obtain the Gagliardo--Nirenberg interpolation inequality \eqref{GN with Boyd indices} via a different method with \cite{LRS2023}.
%Let us replace $\mathbb{T}^d$ with $\mathbb{R}^d$ in Corollary \ref{commutative torus}. 
%Therefore, if the Ces\`{a}ro operator $C$ is bounded on $E(\RR^{d})$ and $F(\RR^{d})$, then the Gagliardo--Nirenberg interpolation inequality \eqref{GN with Boyd indices} holds, i.e.
%\begin{align*}
%\Norm{\nabla^{l}u}_{E(\RR^{d})^{\frac{l}{k}}F(\RR^{d})^{1-\frac{l}{k}}}
%\lesssim\Norm{\nabla^{k}u}_{E(\RR^{d})}^{\frac{l}{k}}\Norm{u}_{F(\RR^{d})}^{1-\frac{l}{k}}
%\quad\mbox{with}\quad u\in W^{k,1}(\RR^{d})\cap F(\RR^{d})
%,~ 0\leq l<k.
%\end{align*}
%Here the assumptions on Boyd indices of symmetric function spaces $E(\RR^{d})$ and $F(\RR^{d})$ are removed, which gives a positive answer to \cite[Question 4.1]{LRS2023}.
\end{remark}

This paper is organised as follows. In Section \ref{Preliminaries}, we introduce noncommutative symmetric spaces, noncommutative torus $\ted$ and the Calder\'{o}n--Lozanovskii space.  In Section \ref{proof}, we present the proof of the Theorem \ref{GN-symmetric}.

\section{Preliminaries}\label{Preliminaries}
\subsection{Noncommutative symmetric spaces}\label{nocommutative symmetric spaces}
Let $\M$ be a von Neumann algebra equipped with a normal faithful normalized trace $\tau$. The set of all $\tau$-measurable operators is denoted to be $L_{0}(\M)$. For every $x\in L_{0}(\M)$, the generalized singular value function of $x$,
\begin{align*}
\mu(x):s\rightarrow\mu(s,x)
\end{align*}
is defined by the formula
\begin{align*}
\mu(s,x)=\inf\set{\Norm{x(1-e)}_{\M}:e\in \proj(\M),~\tau(e)\leq s}.
\end{align*}
The singular value function $\mu(\cdot)$ is sub-additive and $\mu(\cdot,x)$ is non-increasing. For more information about singular value and $\tau$-measurable operators, we refer to \cite{X2007,LSZ2012,LSZM2023}.

\begin{definition}\cite{DPS2024,LSZ2012},
Let $E(\M)$ be a linear subset in $L_{0}(\M)$ equipped with a complete norm $\Norm{\cdot}_{E(\M)}$.  We say that $E(\M)$ is a symmetric operator space if for $x\in E(\M)$ and for every $y\in L_{0}(\M)$ with $\mu(y)\leq\mu(x)$, we have $y\in E(\M)$ and $\Norm{y}_{E(\M)}\leq\Norm{x}_{E(\M)}$.
	
A symmetric function space on $(0,1)$ is the term reserved for a symmetric operator space when $\M=L_{\infty}(0,1)$.
\end{definition}

Recall the construction of noncommutative symmetric (operator) space $E(\M)$. Let $E$ be a symmetric (quasi-)Banach function space on $(0,1)$. Set
$$E(\M)=\set{s\in L_{0}(\M):\mu(x)\in E}.$$ 
We equip $E(\M)$ with a natural (quasi-)norm
$$\Norm{x}_{E(\M)}:=\Norm{\mu(x)}_{E},\quad x\in E(\M).$$
The space $(E(\M),\Norm{\cdot}_{E(\M)})$ is a (quasi-)Banach space with the (quasi-)norm $\Norm{\cdot}_{E(\M)}$ and is called the noncommutative symmetric (quasi-)Banach (operator) space associated with $(\M,\tau)$ corresponding to $(E,\Norm{\cdot}_{E})$ (see \cite{KS2008,S2014}).
The main result of \cite{KS2008,S2014} (see also  \cite[Theorem 3.1.1]{LSZ2012}) shows that the correspondence
$$(E,\Norm{\cdot}_{E})\leftrightarrow(E(\M),\Norm{\cdot}_{E(\M)})$$
is a one-to-one  correspondence between the set of all non commutative symmetric operator spaces in $L_{0}(\M)$ and the set of all symmetric function spaces in $L_{0}(0,1)$.

%A von Neumann algebra $\M$ is atomic if every non-zero projection in $\M$ contains a non-zero minimal projection and atomless if there are no non-zero minimal projections in $\M$. 
%The von Neumann algebras $L_{\infty}((0,1))$ and $L_{\infty}((1,\infty))$ are atomless, while the von Neumann algebra
%$L_{\infty}(\ZZ_{+})$ is atomic (see e.g. \cite[Example 2.1.2]{LSZ2012}). 
%When $\theta=0$, $\ted$ is the commutative von Neumann algebra $L_{\infty}(\TT^{d})$ and it is atomless. {\color{red}Moreover, $\ted$ is atomless.} 

\subsection{Noncommutative torus}\label{noncommutative torus}
The following notation is standard; for more information of noncommutative torus we refer to \cite{CXY2013,LSZM2023,MSX2019,MSZ2019,MSZ2022,SXD2023,XXY2018}.
Fix an integer $d>1$ and a real antisymmetric matrix $\theta=(\theta_{j,k})_{d\times d}$.
The associated $C^{*}$-algebra $C(\ted)$ is the universal $C^{*}$-algebra generated by the $d$ unitary operators $U_{1},U_{2},\cdots,U_{d}$ satisfying the following commutation relation
\begin{align}\label{urelation}
	U_{j}U_{k}=e^{2\pi\ri\theta_{j,k}}U_{k}U_{j},\quad j,k=1,2,\cdots,d.
\end{align}
Let $U=(U_{1},U_{2},\cdots,U_{d})$. For $n=(n_{1},n_{2},\cdots,n_{d})\in\zd$, we denote
$$U^{n}=U_{1}^{n_{1}}U_{2}^{n_{2}}\cdots U_{d}^{n_{d}}.$$
On the space given by the linear span of $\{U^{n}\}_{n\in\ZZ^{d}}$, we define the functional $\tau$ by
\begin{align*}
	\tau(\sum_{n\in\zd}c_{n}U^{n})=c_{\mathbf{0}}
\end{align*}
where $\mathbf{0}= (0,\cdots,0)$. Then, $\tau$ extends to a faithful tracial state on $C(\ted)$. Let $\ted$ be the weak$^{*}$-closure of $C(\ted)$ in the GNS representation of $\tau$. This is our $d$-dimensional non commutative torus. The state $\tau$ extends to a normal faithful tracial state on $\ted$ that will be denoted again by $\tau$. %Recall that the von Neumann algebra $\ted$ is hyperfinite.
Denote $\lp{p}$ $(0<p<\infty)$ the non commutative $L_{p}$-spaces defined by $\tau$ on $\ted$ and denote $\lp{\infty}=\ted$. 

For $x\in\lp{1}$ and $n\in\ZZ^{d}$, we define the Fourier coefficient
$$\wh{x}(n)=\tau(x(U^{n})^{*}).$$
By definition of $\tau$, we see that $\tau(U_{k})=\delta_{k,0}$, and then the standard Hilbert space arguments show that any $x\in\lp{2}$ can be written as an $L_{2}$-convergent series
$\displaystyle x=\sum_{n\in\ZZ^{d}}\wh{x}(n)U^{n}$ with Plancherel identity
$$\Norm{x}_{\lp{2}}^{2}=\sum_{n\in\ZZ^{d}}|\wh{x}(n)|^{2}.$$

The space $C^{\infty}(\ted)$ is defined to be the subset of $x\in C(\ted)$ such that the sequence of Fourier coefficient $\set{\wh{x}(n)}_{n\in\ZZ^{d}}$ has rapid decay of polynomials. We consider $C^{\infty}(\ted)$ as the space of smooth functions on $\ted$. There is also a canonical Fr\'{e}chet topology on $C^{\infty}(\ted)$, and the space $\D'(\ted)$, called the space of distributions on $\ted$, is defined to be the topological dual of $C^{\infty}(\ted)$.

Let $j=1,2,\cdots,d$. The partial differential operator $\p_{j}$ on the linear span of $\{U^{n}\}_{n\in\ZZ^{d}}$ is defined as
\begin{align*}
	\p_{j}(U^{n})=2\pi\ri n_{j}U^{n},\quad n\in\zd.
\end{align*}
Every partial derivation $\p_{j}$ can be viewed as a densely defined closed (unbounded) operator on $\lp{2}$, whose adjoint is equal to $-\p_{j}$.
The gradient operator is denoted by $\nabla=(\p_{1},\cdots,\p_{d})$, which may be considered as a linear operator from $\lp{2}$ to $\lp{2}\otimes\CC^{d}$.
Given $\alpha=(\alpha_{1},\cdots,\alpha_{d})\in\NN^{d}$, we define 
\begin{align*}	\nabla^{\alpha}=\p_{1}^{\alpha_{1}}\cdots\p_{d}^{\alpha_{d}}
	\quad\mbox{and}\quad
	\nabla^{k}=\set{\nabla^{\alpha}}_{\alpha\in\NN^{d},|\alpha|_{1}=k}\quad\mbox{for}\quad k\in\NN.
\end{align*}
Here $|\alpha|_{1}=\alpha_{1}+\cdots\alpha_{d}$. By duality, these derivations transfer to $\D'(\ted)$ as well.
The Sobolev space of order $k\in\NN$ on $\ted$ is defined to be
\begin{align*}
	W^{k,1}(\ted)=\set{x\in\D'(\ted):\nabla^{\alpha} x\in\lp{1}\quad\mbox{with}\quad  |\alpha|_{1}\leq k,~\alpha\in\NN^{d}}.
\end{align*}
More information on the Sobolev spaces, we refer to \cite{XXY2018}.

For our purpose, we use the notation $$D_{j}=-i\p_{j},\quad j=1,2,\cdots,d,$$
which are self-adjoint. Recall  \cite[Chapter 7]{XXY2018} the Fourier multiplier theory on $\ted$. Let $m:\ZZ^{d}\rightarrow\CC$ be a bounded sequence. The Fourier multiplier $m(\nabla)$ is defined in the following way
\begin{align}\label{eitd-multiplier}
m(\nabla)x=\sum_{n\in\ZZ^{d}}m(n)\wh{x}(n)U^{n},\quad x\in\lp{2}.
\end{align} 
The convergence is understood in $L_{2}$-sense. Similarly, we define the operator $m(\nabla_{\td})$ as a Fourier multiplier on $L_{2}(\td)$ with the symbol $m$.

\subsection{Symmetric spaces on $\ted$}

In this paper, we mainly focus on the von Neumann algebra $L_{\infty}(\ted)$. 
Given a symmetric function space $E$ on $(0,1)$, the noncommutative symmetric operator space $E(\ted)$ on $\ted$ is then defined as 
\begin{align*}
E(\ted)=\set{x\in\lo:\mu(x)\in E}\quad\mbox{with}\quad \Norm{x}_{E(\ted)}:=\Norm{\mu(x)}_{E},
\end{align*}
where $\lo$ is the set of all $\tau$-measurable operators on $\ted$.
In the sequel, we always denote $G(\ted)$, $E(\ted)$ and $F(\ted)$ the symmetric operator spaces on $\ted$ related to the symmetric function spaces $G$, $E$ and $F$  on $(0,1)$ separately.

\subsection{The Calder\'{o}n--Lozanovskii space}\label{Calderon-Lozanovskii space}
%Let $L_{0}(0,1)$ be the Lebesgue measurable function space on $(0,1)$. 
Given $0<\epsilon<1$ and two symmetric function spaces $E, F$ on $(0,1)$, the Calder\'{o}n--Lozanovskii space $E^{\epsilon}F^{1-\epsilon}$ is defined as
\begin{align*}
	E^{\epsilon}F^{1-\epsilon}=
	\set{h\in L_{0}(0,1): |h|\leq f^{\epsilon}g^{1-\epsilon}~\mbox{for~some}~ 0\leq f\in E,0\leq g\in F}
\end{align*}
equipped with the norm
\begin{align*}
	\Norm{h}_{E^{\epsilon}F^{1-\epsilon}}
	=\inf\set{\max\set{\Norm{f}_{E},\Norm{g}_{F}}:|h|\leq f^{\epsilon}g^{1-\epsilon},~ 0\leq f\in E,0\leq g\in F}.
\end{align*}
With this norm, $E^{\epsilon}F^{1-\epsilon}$ is a symmetric function space for $0<\epsilon<1$. And it follows from \cite{M1989} that
\begin{align}\label{holder-inequality}
\Norm{|f|^{\epsilon}|g|^{1-\epsilon}}_{E^{\epsilon}F^{1-\epsilon}}
\leq\Norm{f}_{E}^{\epsilon}\Norm{g}_{F}^{1-\epsilon}.
\end{align}
This is the H\"{o}lder's inequality in the Calder\'{o}n--Lozanovskii space. 
%For $1<p<\infty$ and a symmetric function space $E$ on $(0,1)$, define
%\begin{align*}
%E_{p}=\set{f\in L_{0}(0,1):|f|^{p}\in E}
%\end{align*}
%with the (quasi-)norm given by
%\begin{align*}
%\Norm{f}_{E_{p}}=\Norm{|f|^{p}}_{E}^{\frac{1}{p}}.
%\end{align*}
%Observe that $E_{p}=E^{\frac{1}{p}}(L^{\infty})^{1-\frac{1}{p}}$. Therefore, when $E$ is a symmetric function space on $(0,1)$, so are $E_{p}$ $(1<p<\infty)$.
%

\section{Proof of Theorem \ref{GN-symmetric}}\label{proof}
In this section, we prove Theorem \ref{GN-symmetric}.% and Corollary \ref{GN-multipier}.

\begin{lemma}\label{abstract-equality}
Let $A$ be a self-adjoint densely defined closed operator on some Hilbert space. In the strong operator topology, we have
$$iA=e^{iA}-{\rm id}+\int_{0}^{1}\int_{0}^{t} A^2e^{isA}dsdt.$$
\end{lemma}
\begin{proof}
Note that, for $\delta\neq0$,
\begin{align*}
\int_{0}^{t}e^{is\delta}ds=\frac{1}{i\delta}(e^{i\delta t}-1).
\end{align*}
Thus,
\begin{align*}
\int_{0}^{1}\int_{0}^{t}e^{is\delta}dsdt=\frac{1}{i\delta}\int_{0}^{1}(e^{i\delta t}-1)dt
=\frac{1}{(i\delta)^{2}}(e^{i\delta}-1)-\frac{1}{i\delta}.
\end{align*}
Multiplying	both sides of the equality above by $(i\delta)^{2}$ and applying the functional calculus (see e.g.\cite{T2002}), we arrive at the claim.
\end{proof}

\begin{lemma}\label{automorphism}
For $1\leq j\leq d$ and $t\in\RR$, the operator $e^{itD_j}$ is a trace-preserving $\ast$-automorphism of $L_{\infty}(\mathbb{T}^d_{\theta})$.
\end{lemma}
\begin{proof}
Note that $D_{j}(U^{n})=2\pi n_{j}U^{n}$ for any $n\in\ZZ^{d}$. As a Fourier multiplier in \eqref{eitd-multiplier}, we have
\begin{align*}
e^{itD_j}(U^{n})=e^{2\pi it n_{j}}U^{n}.
\end{align*}
Since $e^{2\pi it0}=1$, it follows that $e^{itD_j}$ is trace-preserving.
The above equality and the fact that $\lp{\infty}$ is the weak$^{\ast}$-closure of the linear span of $\{U^{n}\}_{n\in\ZZ^{d}}$ imply that $e^{itD_j}$ is a $\ast$-automorphism  of $L_{\infty}(\mathbb{T}^d_{\theta})$.
\end{proof}

Suppose that $\M_{1}$ and $\M_{2}$ are two von Neumann algebras  equipped with normal faithful normalized traces $\tau_{1}$ and $\tau_{2}$ respectively.
\cite{DPS2024}( see also \cite{LSZ2012,LSZM2023}) If $x\in L_{0}(\M_{1})$ and $y\in L_{0}(\M_{2})$, then $x$ is said to be submajorized by $y$, denoted by $x\dpr y$, if
\begin{align*}
\int_{0}^{t}\mu(s,x)ds\leq\int_{0}^{t}\mu(s,y)ds\quad\mbox{for~all}\quad t>0.
\end{align*} 
The notion of submajorization may, of coursse, in particular be considerd for functions on $(0,1)$, that is, with respect to the commutative von Neumann algebra $L_{\infty}(0,1)$ with Lebesgue integral as trace. If $x\in L_{0}(\M_{1})$ and $y\in L_{0}(\M_{2})$, then it is clear that $y\dpr x$ if and only if $\mu(y)\dpr \mu(x)$.

\begin{lemma}\label{gn key lemma} Let $x\in W^{2,1}(\ted)$. For every $\delta>0,$ we have
$$D_jx\prec\prec \delta^{-1}\mu(x)+\delta\mu(D_j^2x).$$
\end{lemma}
\begin{proof}
Apply Lemma \ref{abstract-equality} to the operator $A=2\delta D_j.$ We have
$$i D_j=(2\delta)^{-1}(e^{2i\delta D_j}-{\rm id})+2\delta\int_0^1\int_0^tD_j^2e^{2is\delta D_j}dsdt.$$
This equality holds pointwise. That is, 
$$i D_jx=(2\delta)^{-1}(e^{2i\delta D_j}x-x)+2\delta\int_0^1\int_0^tD_j^2e^{2is\delta D_j}xdsdt.$$
Using property of submajorization in \cite[Theorem 3.9.9]{DPS2024}, we have
$$D_jx\prec\prec (2\delta)^{-1}(\mu(e^{2i\delta D_j}x)+\mu(x))+2\delta \int_0^1\int_0^t\mu(D_j^2e^{2is\delta D_j}x)dsdt.$$
By Lemma \ref{automorphism}, $e^{2i\delta D_j}$ and $e^{2is\delta D_j}$ are a trace-preserving $\ast$-automorphisms of $L_{\infty}(\mathbb{T}^d_{\theta}).$ In particular \cite[Proposition 3.3 ({\rm ii})]{PS2007}, they preserve the singular value function. Hence, the latter equality is written as
$$D_jx\prec\prec \delta^{-1}\mu(x)+2\delta \int_0^1\int_0^t\mu(D_j^2x)dsdt=\delta^{-1}\mu(x)+\delta\mu(D_j^2x).$$
\end{proof}

\begin{lemma}\label{gn main lemma} Let $x\in W^{2,1}(\ted)$. We have
$$C\mu(D_jx)\leq 2(C\mu(x))^{\frac12}(C\mu(D_j^2x))^{\frac12}.$$
\end{lemma}
\begin{proof} Fix $t\in(0,1).$ By Lemma \ref{gn key lemma}, we have
$$\int_0^t\mu(s,D_jx)ds\leq \delta^{-1}\int_0^t\mu(s,x)ds+\delta \int_0^t\mu(s,D_j^2x)ds$$
for every $\delta>0.$ Since $x\in W^{2,1}(\ted)$, the above two integrals on the right hand side are finite. If $\displaystyle\int_0^t\mu(s,D_j^2x)ds$ is zero, letting $\delta\rightarrow+\infty$ yields the desired result. If $\displaystyle\int_0^t\mu(s,D_j^2x)ds\neq0$, setting
$$\delta=\Big(\frac{\int_0^t\mu(s,x)ds}{\int_0^t\mu(s,D_j^2x)ds}\Big)^{\frac12},$$
we obtain
$$\int_0^t\mu(s,D_jx)ds\leq 2\Big(\int_0^t\mu(s,x)ds\Big)^{\frac12}\Big( \int_0^t\mu(s,D_j^2x)ds\Big)^{\frac12}.$$
Since $t\in(0,1)$ is arbitrary, the assertion follows.
\end{proof}

The following sequence is proposed in the main proof for Theorem \ref{GN-symmetric}. For our purpose, let us firstly investigate the iterative relations of the first term and any other term. 
\begin{lemma}\label{geometric mean lemma} Suppose $(a_l)_{l\geq 0}$ are positive numbers such that
$$a_l\leq 2a_{l-1}^{\frac12}a_{l+1}^{\frac12},\quad l\in\mathbb{N}.$$
It follows that
$$a_l\leq 2^{3\cdot 2^{k-2}-2}a_0^{1-\frac{l}{k}}a_k^{\frac{l}{k}},\quad 0\leq l\leq k,\quad k\in\NN.$$
\end{lemma}
\begin{proof} Denote for brevity $c_k=2^{3\cdot 2^{k-2}-2},$ $k\geq 2.$

We prove the assertion by induction on $k.$ Base of induction (i.e., the case $k=2$) is obvious. Let us establish the step of the induction. Suppose the assertion holds for $k$ and let us establish it for $k+1.$

We have
$$a_k\leq 2a_{k-1}^{\frac12}a_{k+1}^{\frac12}.$$
By inductive assumption, we have
$$a_{k-1}\leq c_ka_0^{\frac1k}a_k^{1-\frac1k}.$$
Thus,
$$a_k\leq 2c_k^{\frac12}a_0^{\frac1{2k}}a_k^{\frac12-\frac1{2k}}a_{k+1}^{\frac12}.$$
Either $a_k=0$ or, if $a_k$ is strictly positive, then
$$a_k^{\frac12+\frac1{2k}}\leq 2c_k^{\frac12}a_0^{\frac1{2k}}a_{k+1}^{\frac12}.$$
Taking both sides to the power $\frac{2k}{k+1},$ we obtain
$$a_k\leq (4c_k)^{\frac{k}{k+1}}a_0^{\frac1{k+1}}a_{k+1}^{1-\frac1{k+1}}.$$

By inductive assumption, we have
$$a_l\leq c_ka_0^{1-\frac{l}{k}}a_k^{\frac{l}{k}},\quad 0\leq l\leq k.$$
Thus,
$$a_l\leq c_ka_0^{1-\frac{l}{k}}\cdot (4c_k)^{\frac{k}{k+1}\cdot\frac{l}{k}}a_0^{\frac1{k+1}\cdot\frac{l}{k}}a_{k+1}^{(1-\frac1{k+1})\cdot \frac{l}{k}}=c_k\cdot (4c_k)^{\frac{l}{k+1}}\cdot a_0^{1-\frac{l}{k+1}}a_{k+1}^{\frac{l}{k+1}},\quad 0\leq l\leq k.$$
Since $4c_k\geq 1,$ it follows that
$$a_l\leq 4c_k^2a_0^{1-\frac{l}{k+1}}a_{k+1}^{\frac{l}{k+1}}=c_{k+1}a_0^{1-\frac{l}{k+1}}a_{k+1}^{\frac{l}{k+1}},\quad 0\leq l\leq k+1.$$
This yields the step of the induction and, hence, completes the proof.
\end{proof}

\begin{lemma}\label{second geometric mean lemma} Let $a:\mathbb{Z}^d_+\to\mathbb{R}_+$ be such that
$$a(\beta+e_j)\leq 2a(\beta)^{\frac12}a(\beta+2e_j)^{\frac12},\quad \beta\in\mathbb{Z}^d_+,\quad 1\leq j\leq d.$$
For every $\alpha\in\mathbb{Z}^d$ and for every $k\geq|\alpha|_1,$ we have
$$a(\alpha)\leq 2^{3\cdot 2^{k-2}-2} a(0)^{1-\frac{|\alpha|_1}{k}}\cdot\big(\max_{|\beta|_1=k}a(\beta)\big)^{\frac{|\alpha|_1}{k}}.$$
\end{lemma}
\begin{proof} 
Set
$$a_m=\max_{|\beta|_1=m}a(\beta),\quad m\in\mathbb{Z}_+.$$ 

Denote $|\alpha|_{1}=l$. For $l\in\mathbb{N},$ find $\beta\in\mathbb{Z}^d_+$ such that $|\beta|_1=l$ and such that $a_l=a(\beta).$ Fix $1\leq j\leq d$, we have
$$a(\alpha)\leq a_l=a(\beta)\leq 2a(\beta-e_j)^{\frac12}a(\beta+e_j)^{\frac12}\leq 2a_{l-1}^{\frac12}a_{l+1}^{\frac12},$$
where we use the assumption in the second inequality.
The assertion follows now from Lemma \ref{geometric mean lemma}.
\end{proof}

\begin{lemma}\label{gn post-main lemma} Let $x\in W^{k,1}(\ted)$. For every $\alpha\in\mathbb{Z}^d_+$ and for every $k\geq|\alpha|_1,$ we have
$$C\mu(\nabla^{\alpha}x)\leq 2^{3\cdot 2^{k-2}-2}(C\mu(x))^{1-\frac{|\alpha|_1}{k}}(\max_{|\beta|_1=k}C\mu(\nabla^{\beta}x))^{\frac{|\alpha|_1}{k}}.$$
\end{lemma}
\begin{proof} Applying Lemma \ref{gn main lemma} with $\nabla^{\alpha}x$ instead of $x,$ we obtain
$$C\mu(\nabla^{\alpha+e_j}x)\leq 2(C\mu(\nabla^{\alpha}x))^{\frac12}(C\mu(\nabla^{\alpha+2e_j}x))^{\frac12},\quad \alpha\in\mathbb{Z}^d_+,\quad 1\leq j\leq d.$$
Fix $t\in (0,1)$ and set
$$a(\alpha)=(C\mu(\nabla^{\alpha}x))(t),\quad \alpha\in\mathbb{Z}^d_+.$$
We have
$$a(\alpha+e_{j})\leq2a(\alpha)^{\frac12}a(\alpha+2e_{j})^{\frac12}.$$
Applying Lemma \ref{second geometric mean lemma}, we obtain
$$(C\mu(\nabla^{\alpha}x))(t)\leq 2^{3\cdot 2^{k-2}-2} \big((C\mu(x))(t)\big)^{1-\frac{|\alpha|_1}{k}}\cdot\big(\max_{|\beta|_1=k}(C\mu(\nabla^{\beta}x))(t)\big)^{\frac{|\alpha|_1}{k}}.$$
\end{proof}

\begin{proof}[Proof of Theorem \ref{GN-symmetric}] Let $x\in W^{k,1}(\ted)$. By Lemma \ref{gn post-main lemma}, we have
$$\mu(\nabla^{\alpha}x)\leq 2^{3\cdot 2^{k-2}-2}(C\mu(x))^{1-\frac{|\alpha|_1}{k}}(\max_{|\beta|_1=k}C\mu(\nabla^{\beta}x))^{\frac{|\alpha|_1}{k}}.$$
Thus, combining the above estimate and H\"{o}lder's inequality  \eqref{holder-inequality}, we obtain
\begin{align*}
\|\nabla^{\alpha}x\|_{G(\ted)}
&\leq 2^{3\cdot 2^{k-2}-2}\|(C\mu(x))^{1-\frac{|\alpha|_1}{k}}(\max_{|\beta|_1=k}C\mu(\nabla^{\beta}x))^{\frac{|\alpha|_1}{k}}\|_{E^{1-\frac{|\alpha|_1}{k}}F^{\frac{|\alpha|_1}{k}}}\\
&\leq 2^{3\cdot 2^{k-2}-2}\|C\mu(x)\|_E^{1-\frac{|\alpha|_1}{k}}\|\max_{|\beta|_1=k}C\mu(\nabla^{\beta}x)\|_F^{\frac{|\alpha|_1}{k}}\\
&\leq 2^{3\cdot 2^{k-2}-2}\|C\mu(x)\|_E^{1-\frac{|\alpha|_1}{k}}\big(\sum_{|\beta|_1=k}\|C\mu(\nabla^{\beta}x)\|_F\big)^{\frac{|\alpha|_1}{k}}\\
&\leq 2^{3\cdot 2^{k-2}-2}\|C\|_{E\to E}^{1-\frac{|\alpha|_1}{k}}\|C\|_{F\to F}^{\frac{|\alpha|_1}{k}}\|x\|_{E(\ted)}^{1-\frac{|\alpha|_1}{k}}\big(\sum_{|\beta|_1=k}\|\nabla^{\beta}x\|_{F(\ted)}\big)^{\frac{|\alpha|_1}{k}}.
\end{align*}
This completes the proof.
\end{proof}

{\bf Acknowledgements:} F. Sukochev  and D. Zanin are supported by Australian Research Council (ARC) DP230100434. F. L. Yang is supported by National Natural Science Foundation of China No. 12371138 and he is grateful for the help of professor Jinghao Huang and professor Xiao Xiong.

\end{document}